\documentclass[12pt,a4paper]{article}

\usepackage{amsmath,amssymb,amsfonts,mathtools,amsthm}
\usepackage{tikz}
\usepackage{algorithm}
\usepackage{algpseudocode}
\usepackage{graphicx,epsfig,color,scalerel,subfigure,booktabs}
\usepackage[super]{nth}
\usepackage{xparse}
\usepackage{siunitx}
\usepackage{afterpage}
\usepackage{enumitem}
\usepackage{comment}

\newtheorem{lemma}{Lemma}
\newtheorem{theorem}{Theorem}

\definecolor{lightblue}{rgb}{0.22,0.45,0.70}
\usepackage[colorlinks=true,breaklinks=true,linkcolor=lightblue,citecolor=lightblue]{hyperref}

\numberwithin{equation}{section}
\numberwithin{figure}{section}
\numberwithin{table}{section}

\newcommand{\bn}{\boldsymbol{n}}

\newcommand{\bu}{\boldsymbol{u}}
\newcommand{\bv}{\boldsymbol{v}}
\newcommand{\bw}{\boldsymbol{w}}

\newcommand{\fb}{\boldsymbol{f}}

\newcommand{\bx}{\boldsymbol{x}}
\newcommand{\by}{\boldsymbol{y}}

\newcommand\beps{\boldsymbol{\varepsilon}}

\newcommand\bzeta{\boldsymbol{\zeta}}
\newcommand\bzero{\boldsymbol{0}}
\newcommand\bxi{\boldsymbol{\xi}}

\newcommand\bPi{\boldsymbol{\Pi}}

\newcommand{\bH}{\mathbf{H}}

\newcommand{\bL}{\mathbf{L}}

\newcommand{\bV}{\mathbf{V}}
\newcommand{\bQ}{\mathbf{Q}}

\newcommand{\bbM}{\mathbb{M}}
\newcommand{\bbI}{\mathbb{I}}
\newcommand{\bbR}{\mathbb{R}}

\newcommand\bdiv{\mathop{\mathbf{div}}\nolimits}
\newcommand\vdiv{\mathop{\mathrm{div}}\nolimits}

\newcommand\bcurl{\mathop{\mathbf{curl}}\nolimits}

\newcommand\tr{\mathop{\mathrm{tr}}\nolimits}

\DeclarePairedDelimiter\norm{\lVert}{\rVert}

\def\XXint#1#2#3{{\setbox0=\hbox{$#1{#2#3}{\int}$ }
		\vcenter{\hbox{$#2#3$ }}\kern-.6\wd0}}

\title{Robust virtual element methods for 3D stress-assisted diffusion problems}

\author{Andr\'es E. Rubiano\thanks{School of Mathematics, Monash University, 9 Rainforest Walk, 3800 VIC, Australia.  \protect\url{andres.rubianomartinez@monash.edu},
\protect\url{https://orcid.org/0000-0002-5557-4963}
}}

\date{\today}

\usepackage[square,sort,comma,numbers]{natbib}
\begin{document}

\maketitle

\begin{abstract}
This paper presents an initial exploration of stress-assisted diffusion problems in three dimensions within the framework of the virtual element method (VEM). Hilbert spaces enriched with parameter-weighted norms, the extended Babu$\check{\text{s}}$ka–Brezzi–Braess theory for perturbed saddle-point
problems, and Banach fixed-point theory play a crucial role in performing a robust analysis of the fully coupled non-linear system. The proposed virtual element formulations are provided with appropriate projection, interpolation, and stabilisation operators that ensures the well-posedness of the discrete problem. Numerical simulations are conducted to show the accuracy, performance, and applicability of the method.

\end{abstract}

\section{Introduction}
\paragraph{Scope.} Stress-assisted diffusion appears broadly in applications like Lithium-ion battery cells, silicon rubber and hydrogen diffusion, polymer-based coatings, semiconductor fabrication, oxidation of silicon nanostructures, and enhancing of conductivity properties in soft living tissue. In these processes, the stress of the material involved affects the patterns of diffusion. This paper addresses one of the simplest models for such an interaction, incorporating stress effects inside the diffusion coefficient, and recovering Fick's law in the absence of these phenomena (see, e.g., \cite{cherubini17}). Moreover, we extend the analysis done in \cite{khot2024} pointing out the definition of VEM spaces in 3D and the polynomial, interpolation, and stabilisation operators used in this framework.

\paragraph{Outline.} The content of the paper is organised as follows: The remainder of this section presents the strong form of the coupled fully coupled stress-assisted diffusion problem. Section~\ref{sec:weak-formulation} is devoted to presenting the weak formulation and well-posedness analysis using perturbed saddle-point theory and fixed-point arguments. In Section~\ref{sec:vem} we define the discrete problem, show the existence and uniqueness of discrete solution, and outline the a priori error analysis. Finally, in Section~\ref{sec:results} we provide numerical results that illustrate the theoretical results presented, and we present as an application of the model the lithiation of a perforated cylindrical anode particle. 


\paragraph{Model statement.}
Let us consider the following nonlinear formulation for stress-assisted diffusion of a solute that interacts with an elastic material
\begin{subequations}\label{eq:mixed-formulation}
\begin{align}
-\bdiv(2\mu \beps(\bu) -p\bbI) = \fb &\quad \text{in $\Omega$},\quad \bu=\mathbf{0} \quad \text{on $\Gamma_D$},\label{eq:linear-momentum}\\
p = -\lambda \vdiv\bu\, +\ell(\varphi) &\quad \text{in $\Omega$},\quad (2\mu \beps(\bu) -p\bbI)\bn = \bzero \quad \text{on $\Gamma_N$},\label{eq:hooke-law}\\
\bzeta = \bbM(\beps(\bu),p) \nabla \varphi &\quad \text{in $\Omega$},\quad \varphi=\varphi_D \quad \text{on $\Gamma_D$},\label{eq:diffusive-flux}\\
\theta \varphi 
- \vdiv(\bzeta) = g &\quad \text{in $\Omega$},\quad \bzeta \cdot \bn = 0 \quad \text{on $\Gamma_N$}\label{eq:reaction-diffusion}.
\end{align}
\end{subequations}
Equations \eqref{eq:linear-momentum} and \eqref{eq:hooke-law} constitute a Herrmann-like mixed formulation for linear elasticity. The coupling between the displacement $\bu$ of the material, the pressure $p$, and the solute concentration $\varphi$ follows Hooke's law to express the Cauchy stress tensor as $2\mu \beps(\bu) -p\bbI$, $\mu$ and $\lambda$ are the Lam\'e parameters, $\ell$ is the active stress coefficient, and $\fb$ is a vector of external body loads. 

On the other hand, equations \eqref{eq:reaction-diffusion} and \eqref{eq:diffusive-flux} correspond to the reaction-diffusion equation written in mixed form where $\bzeta$ is the diffusive flux, $g$ is a given net volumetric source of solute, $\theta$ is a positive model parameter, and $\bbM(\beps(\bu),p)$ is the stress-assisted diffusion coefficient. 

We adopt mixed loading boundary conditions for the coupled problem: the structure is clamped and has a given concentration on $\Gamma_D$,  where the boundary subset $\Gamma_D\subset  \partial \Omega$ is of positive surface measure. In addition, a given traction and zero solute flux are prescribed on $\Gamma_N:= \partial\Omega \setminus \Gamma_D$. 

\paragraph{The nonlinear terms.} We assume that $\bbM(\cdot,\cdot)$ is symmetric, positive semi-definite and uniformly bounded in $\mathbb{L}^\infty(\Omega)$, likewise for $\bbM^{-1}(\cdot,\cdot)$. More explicitly, for all $\bw\in \bH^1(\Omega), r\in \text{L}^2(\Omega) \text{ and }  \bx,\by \in \mathbb{R}^3$  there exists $M\in \mathbb{R}$ such that $0<M^{-1} \leq M$ with 
$M^{-1} \bx\cdot \bx \leq \bx\cdot[\bbM^{-1}(\beps(\bw),r)\bx]$, and $\by\cdot[\bbM^{-1}(\beps(\bw),r)\bx]  \leq M \bx\cdot\by$. In addition, we assume that $\ell: \text{L}^2(\Omega)\to \text{L}^2(\Omega)$ and satisfies $\norm{\ell(\vartheta)}_{0,\Omega} \lesssim \norm{\vartheta}_{0,\Omega}$ for all $\vartheta \in \text{L}^2(\Omega)$. Moreover, we assume that $\bbM^{-1}(\cdot,\cdot)$ and $\ell(\cdot)$ are Lipschitz continuous with Lipschitz constants $L_{\bbM}$ and $L_{\ell}$. Examples of these terms can be found in   \cite{grigoreva19} for the stress-assisted diffusion and \cite{murray2003mathematical} for the active stress.
\section{Weak formulation}\label{sec:weak-formulation}
In view of the boundary conditions, we define the  Hilbert spaces 
\begin{gather*}
\bH^1_D(\Omega):=\{\bv \in \bH^1(\Omega): \bv = \bzero\, \text{on }\Gamma_D\}, \\ 
\bH_N(\vdiv,\Omega):=\{\bxi \in \bH(\vdiv,\Omega): \bxi\cdot\bn = 0 \, \text{on }\Gamma_N\},   
\end{gather*}
with the boundary assignment understood in the sense of traces, and consider the following weak formulation: for given $\fb\in \bL^2(\Omega)$, $g\in \text{L}^2(\Omega)$, 
and $\varphi_D\in H^{1/2}(\Gamma_D)$, find 
$(\bu,p,\bzeta,\varphi) \in \bH_D^1(\Omega)\times \text{L}^2(\Omega) \times \bH_N(\vdiv,\Omega) \times \text{L}^2(\Omega)$ such that 
\begin{subequations}\label{eq:weak}
\begin{gather}
 2\mu \int_\Omega \beps(\bu):\beps(\bv) - \int_\Omega p\vdiv\bv  = \int_\Omega \fb\cdot\bv, \; \forall \bv \in \bH_D^1(\Omega) \\
- \int_\Omega q\vdiv\bu - \lambda^{-1} \int_\Omega p q  = \lambda^{-1}\int_\Omega  \ell(\varphi)q,\; \forall q \in \text{L}^2(\Omega)\\
\int_\Omega \bbM(\beps(\bu),p)^{-1} \bzeta \cdot \bxi + \int_\Omega \varphi \vdiv \bxi  = \langle \varphi_D, \bxi\cdot\bn\rangle_{\Gamma_D},\; \forall \bxi \in \bH_N(\vdiv,\Omega)\\ 
\int_\Omega \psi \vdiv\bzeta - \theta \int_\Omega \varphi\psi  
= - \int_\Omega g\psi, \; \forall \psi \in \text{L}^2(\Omega).
\end{gather}
\end{subequations}
\paragraph{Unique solvability with parameter-weighted norms.} Let us adopt the following notation for the functional spaces for displacement and total volumetric stress $\bV_1 := \bH_D^1(\Omega)$ and $Q_1 = Q_{b_1} := \text{L}^2(\Omega)$, equipped with the scaled norms and semi-norms given by
\begin{gather*}
    \norm{\bu}_{\bV_1}^2 := 2\mu\norm{\beps(\bu)}_{0,\Omega}^2, \quad 
    \norm{p}_{Q_1}^2 := \left((2\mu)^{-1} + \lambda^{-1}\right)\norm{p}_{0,\Omega}^2,\\
    |\bv|_{1,\bV_1}^2 := 2\mu|\bv|_{1,\Omega}^2, \quad |q|_{1,Q_{b_1}}^2 := (2\mu)^{-1}|q|_{1,\Omega}^2.
\end{gather*}
On the other hand, let us denote the functional spaces for diffusive flux and concentration as $\bV_2 = \bH_N(\vdiv,\Omega)$ and $Q_2 = Q_{b_2} := \text{L}^2(\Omega)$, furnished with the following norms and semi-norms as $\norm{\bzeta}_{\bbM,\Omega}^2 := \int_\Omega \bbM(\beps(\bu),p)^{-1} \bzeta \cdot \bzeta$,
\begin{gather*}
    \norm{\bzeta}_{\bV_2}^2 := \norm{\bzeta}_{\bbM,\Omega}^2 + M\norm{\vdiv \bzeta}_{0,\Omega}^2,
    \quad
    \norm{\varphi}_{Q_2}^2 := \left(M^{-1} + \theta \right)\norm{\varphi}_{0,\Omega}^2,\\
    |\bxi|_{1,\bV_2}^2 := M|\bxi|_{1,\Omega}^2, \quad |\psi|_{1,Q_{b_2}}^2 := M^{-1}|\psi|_{1,\Omega}^2.
\end{gather*}
The proposed spaces in conjunction with the extended Babuška--Brezzi--Braess theory for perturbed saddle-point problems and a fixed-point argument show the robust unique solvability of our weak formulation presented in Section \ref{sec:weak-formulation}. We finish this section by stating the continuous dependence on data (see \cite[Section 3]{khot2024} for details).

\begin{theorem}\label{well-posedness}
Let $W =\{ w \in Q_2 \colon \norm{w}_{Q_2} \leq C_2 \sqrt{M}  \norm{\varphi_D}_{1/2,\Gamma_D} + \norm{g}_{0,\Omega}) \}$. Under the assumptions over the non-linear terms, suppose further that $1\leq \lambda$, $1\leq\mu$, $\theta \leq M^{-1}$, and $C_1 L_\ell \sqrt{2\mu} M^{2}C_2^2 L_{\bbM}(\norm{\varphi_D}_{1/2,\Gamma_D} + \norm{g}_{0,\Omega}) < 1.$
Then, for $\varphi \in W$ there is an unique solution $(\bu,p,\bzeta,\varphi)\in \bV_1\times Q_1 \times \bV_2 \times Q_2$ of \eqref{eq:weak} such that 
\begin{subequations}
    \begin{align}
        \norm{(\bu,p)}_{\bV_1\times Q_1} &\leq C_1 \left( \norm{F_1}_{\bV'_1} + \norm{G^\varphi_1}_{Q'_1} \right),\label{continuos-dependence-elasticity}\\
        \norm{(\bzeta,\varphi)}_{\bV_2\times Q_2} &\leq C_2 \left( \norm{F_2}_{\bV'_2} +\norm{G_2}_{Q'_2} \right)\label{continuos-dependence-diffusion},
    \end{align}
\end{subequations}
where the constants $C_1$ and $C_2$ do not depend on the physical parameters. 
\end{theorem}


\section{Virtual element  discretisation}\label{sec:vem}
\paragraph{Mesh assumptions.} Let $\mathcal{T}^h$ be a decomposition of $\Omega$ into polyhedral elements $P$ with  diameter $h_P$, let $\mathcal{F}^h$ be the set of faces $f$ with length $h_f$ and $\mathcal{E}^h$ be the set of edges $e$ of $\mathcal{T}^h$ with length $h_e$. The following mesh assumptions are considered throughout this paper. We assume that there exists a universal constant $\rho>0$ such that
\begin{enumerate}[label={(\bfseries M\arabic*)}]
    \item \label{M1} Each polyhedral element $P$ is star-shaped with respect to a ball of radius $\geq$ $\rho h_P$,
    \item \label{M2} Every face $f$ of $P$ is of diameter $h_f$ and is star-shaped with respect to a disk of radius $\geq$ $\rho h_P$,
    \item \label{M3} Every edge $e$ of $P$ has length $\geq$ $\rho h_P$. 
\end{enumerate}

\paragraph{Polynomial spaces.} Given  an integer $k\geq 0$ the space of polynomials of degree $\leq k$ on $P$ is denoted by $\mathcal{P}_k(P)$ (resp. for faces $f$). The space of the gradients of polynomials of grade $\leq k+1$ on $P$ is denoted as $\mathcal{G}_k(P):= \nabla (\mathcal{P}_{k+1}(P))$ with standard notation $\mathcal{P}_{-1}(P)=\{0\}$ for $k=-1$. The space $\mathcal{G}_k^\oplus(P)$ denotes the complement of the space $\mathcal{G}_k(P)$ in the vector polynomial space $(\mathcal{P}_k(P))^3$,  that is, $(\mathcal{P}_k(P))^3 = \mathcal{G}_k(P) \oplus  \mathcal{G}_k^\oplus(P)$. In particular, following \cite{veiga19}, we set $\mathcal{G}_k^\oplus(P):= \bx\wedge(\mathcal{P}_{k-1}(P))^3$ with $\bx = (x_1,x_2,x_3)^{\tt t}$. Likewise, the space that defines the rotational of polynomials with degree $\leq k+1$ is denoted as $\mathcal{R}_{k}(P):= \bcurl(\mathcal{P}_{k+1}(P))$ where the associated complement space $\mathcal{R}_k^\oplus(P)$ fulfills the property $(\mathcal{P}_k(P))^3 = \mathcal{R}_k(P) \oplus  \mathcal{R}_k^\oplus(P)$ with $\mathcal{R}_k^\oplus(P) = \bx \mathcal{P}_{k-1}(P)$.

Let $\bx_P = (x_{1,P},x_{2,P},x_{3,P})^{\tt t}$ denote the barycentre of $P$ and let $\mathcal{M}_k(P)$ be the set of scaled monomials
\[\mathcal{M}_k(P):=\left\{ \left( \frac{\bx-\bx_P}{h_P} \right)^{\boldsymbol{\alpha}}, 0\leq |\boldsymbol{\alpha}|\leq k \right\},\]
where $\boldsymbol{\alpha}=(\alpha_1,\alpha_2,\alpha_3)^{\tt t}$ is a non-negative multi-index with $|\boldsymbol{\alpha}|=\alpha_1+\alpha_2+\alpha_3$ and $\bx^{\boldsymbol{\alpha}}=x_1^{\alpha_1}x_2^{\alpha_2}x_3^{\alpha_3}$. In particular, we can take the basis of $\mathcal{G}_{k}(P)$ and $\mathcal{G}_{k}^\oplus(P)$ as $\mathcal{M}_{k}^{\nabla}(P):=\nabla\mathcal{M}_{k+1}(P)\setminus\{\mathbf{0}\}$ and $\mathcal{M}_{k}^{\oplus}(P):=\bx_P \wedge (\mathcal{M}_{k-1}(P))^3$, where $(\mathcal{M}_k(P))^3 = \mathcal{M}_k^{\nabla}(P) \oplus  \mathcal{M}_k^{\oplus}(P)$ holds. In addition, we introduce the notation $\mathcal{M}_{k \setminus k-1}(P):= \mathcal{M}_{k}(P) \setminus\mathcal{M}_{k-1}(P)$.

\paragraph{Discrete formulation for the elasticity problem.} The definition of the enhanced 3D VE space follows the approach for Stokes-like problems given in \cite{beirao2020stokes}. Given $k_1\geq 2$. First, we define an extended enhance local VE space as
\begin{align*}
    \bV^{h,k_1}_1(P) := & \{
            \bv_h\in \bH^1(P) \colon \bv_h|_{\partial P}\in (\widetilde{\mathcal{B}}^{h,k_1}_1(\partial P))^3,\;\vdiv \bv_h\in\mathcal{P}_{k_1-1}(P), \\
            & \; -2\mu \bdiv\beps(\bv_h) -\nabla s\in \mathcal{G}_{k_1}^\oplus(P) \text{ for some } s\in \text{L}^2_0(P),\\
            & \; \int_P (\bv_h-\bPi_1^{\beps,k_1}\bv_h)\cdot \mathbf{m}_{k_1}^\oplus = 0, \forall \mathbf{m}_{k_1}^\oplus\in \mathcal{M}_{k_1\setminus k_1-2}^\oplus(P)\},
\end{align*}
where the boundary space of VE functions along the boundary $\partial P$ of $P$, is defined as follows
$$\widetilde{\mathcal{B}}^{h,k_1}_1(\partial P):=\{v_h\in C^0(\partial P)\colon v_h|_f\in \widetilde{\mathcal{B}}^{h,k_1}_1(f), \, \forall f\subset\partial P\},$$
and for each face $f\in\partial P$, the enhanced VE space $\widetilde{\mathcal{B}}^{h,k_1}_1(f)$ locally solves the Poisson equation with Dirichlet boundary conditions and is defined by 
\begin{align*}
    \widetilde{\mathcal{B}}^{h,k_1}_1(f):= &\{
        v_h\in H^1(f) \colon v_h|_{\partial f}\in C^0(\partial f),\; v|_e\in\mathcal{P}_{k_1}(e), \; \Delta_f v_h\in\mathcal{P}_{k_1+1}(f), \\ & \; \hspace{-0.2cm} \int_f (v_h-\Pi_1^{\beps,k_1,f}v_h) m_{k_1+1} = 0, \,  \forall e\subset\partial f, \, \forall m_{k_1+1}\in\mathcal{M}_{k_1+1 \setminus k_1-2}(f)\},
\end{align*}
where $\Delta_f$ denotes the tangential differential operator on $f$, $\Pi_1^{\beps,k_1,f}$ is the restriction to the face $f$ of the energy projection operator for scalar functions defined in \eqref{energy-proj}. The global discrete spaces are set as
\begin{gather*}
\bV^{h,k_1}_1 := \{\bv_h\in \bV_1: \bv_h|_P\in\bV^{h,k_1}_1(P), \, \forall P\in\mathcal{T}^h\}, \\ Q_1^{h,k_1} := \{ q_h\in Q_1 \colon q_h|_P\in \mathcal{P}_{k_1-1}(P), \, \forall P\in \mathcal{T}^h \}.
\end{gather*}
\paragraph{Discrete formulation for the reaction-diffusion problem.} The construction of $\bH(\vdiv,\Omega)$ conforming 3D VE space naturally follows the same approach as its 2D counterpart. For further details, we refer the reader to \cite{veiga-Hdiv}. In this case, we differ from the 2D version by setting $k_2\geq 1$. The discrete VE space locally solve a $\nabla(\vdiv)-\bcurl$ problem as follows
\begin{align*}
    \bV_2^{h,k_2}(P) := &\{ \bxi_h \in \bH(\vdiv,P)\cap \bH(\bcurl,P) \colon \bxi_h \cdot \bn_P^f|_f \in \mathcal{P}_{k_2}(f), \, \forall f\in \partial P, \\
    & \; \nabla (\vdiv \bxi_h) \in \mathcal{G}_{k_2-2}(P), \bcurl \bxi_h \in \mathcal{R}_{k_2-1}(P)\}.
\end{align*}
Then, the discrete global spaces are defined by
\begin{gather*}
    \bV_2^{h,k_2} := \{ \bxi_h \in \bV_2 \colon \bxi_h|_P \in \bV_2^{h,k_2}(P), \, \forall P\in \mathcal{T}^h \}, \\ Q_2^{h,k_2} := \{ \psi_h\in Q_2 \colon \psi_h|_P\in \mathcal{P}_{k_2-1}(P), \, \forall P\in \mathcal{T}^h \}.
\end{gather*}
\paragraph{Polynomial projections.} The energy projection operator given by $\bPi_{1}^{\beps,k_1}: \bH^1(P)\rightarrow (\mathcal{P}_{k_1}(P))^3$, and the $\text{L}^2$-projections operators defined as $\bPi_{j}^{0,k_j}: \bL^2(P)\rightarrow (\mathcal{P}_{k_j}(P))^3$ with $j=1,2$ fulfill that
\begin{align}
    &\bullet \begin{dcases}\int_P \beps(\bv_h-\bPi_{1}^{\beps,k_1}\bv_h) \colon \beps(\mathbf{m}_{k_1}) = 0, \quad \forall \mathbf{m}_{k_1} \in (\mathcal{M}_{k_1}(P))^3,\\
    \int_{\partial P} (\bPi_{1}^{\beps,k_1} \bv - \bv) \cdot \mathbf{m}_{\text{RBM}} =  0, \quad \forall\mathbf{m}_{\text{RBM}}\in \text{RBM}(P),\end{dcases}\label{energy-proj}\\
    &\bullet \int_P (\bv_h-\bPi_{j}^{0,k_j}\bv_h)\cdot \mathbf{m}_{k_{j}} = 0, \quad \forall \mathbf{m}_{k_{j}} \in (\mathcal{M}_{k_{j}}(P))^3.,
\end{align}
where $\text{RBM}(P)$ denotes the set of scaled rigid body motions. We recall that the definition for scalar functions is analogous with the usual notation $\Pi_1^{\beps,k_1}$ and $\Pi_j^{0,k_j}$ for $j=1,2$. Regarding computability, we refer the reader to \cite[Proposition 5.1]{beirao2020stokes} and \cite[Theorem 3.2]{veiga-Hdiv}.
\paragraph{Interpolation operators.} Let $1/2<s_1\leq k_1$, $0\leq s_2 \leq k_2$. The Fortin--like (see \cite{daveiga2022stability,dassi2022bend}) interpolation operators
$\bPi_1^{F,k_1}:\mathbf{H}^{s_1+1}(P)\rightarrow \bV_1^{h,k_1}(P)$, and $\bPi_2^{F,k_2}:\mathbf{H}^{s_2+1}(P)\rightarrow \bV_2^{h,k_2}(P)$ are defined through the Degrees of freedom as
\begin{align*}
    \bullet \, \mathrm{DoF}_{j}(\bv-\bPi_1^{F,k_j}\bv)=0, \quad \forall \bv\in\mathbf{H}^{s_j+1}(P),\; j=1,\dots, \dim(\bV_j^{h,k_j}(P)),
\end{align*}
where the operator $\mathrm{DoF}_j$ indicates the application of the $j$-th Degree of Freedom (see \cite{beirao2020stokes,veiga-Hdiv}). Note that the associated commutative property is given as $\vdiv \bPi_j^{F,k_j}(\cdot) =  \Pi_j^{0,k_j-1} \vdiv (\cdot).$
\paragraph{Approximation and interpolation estimates.} This paragraph collects all the results involving approximation and interpolation estimates needed for the discrete formulation. We recall that arise as consequence of classical estimates and for the robust version of the estimates we refer to \cite[Section 4]{khot2024} pointing out that the extension to the 3D estimates follow similarly.
\begin{lemma} \label{approximation-estimates}
    For any $\bv \in (\bH^{s_1+1}(P)\cap \bV_1(P),|\cdot|_{1,\bV_1(P)})$, $q\in (H^{s_1+1}(P)\cap Q_{b_1}(P),|\cdot|_{1,Q_{b_1}(P)})$, $\bxi \in (\bH^{s_2+1}(P)\cap \bV_2(P),|\cdot|_{1,\bV_2(P)})$ and $\psi \in (H^{s_2+1}(P)\cap Q_{b_2}(P),|\cdot|_{1,Q_{b_2}(P)})$, the polynomial projections $\bPi_1^{\beps,k_1}\bv$, $\Pi_1^{0,k_1}q$, $\bPi_2^{0,k_2}\bxi$ and $\Pi_2^{0,k_2}\psi$ satisfy the following estimates
    \begin{gather*}
    \norm{\bv-\bPi_1^{\beps,k_1}\bv}_{\bV_1(P)}\lesssim h_P^{s_1}|\bv|_{\bV_1(P)},\quad \norm{q-\Pi_1^{0,k_1}q}_{Q_{b_1}(P)}\lesssim h_P^{s_1+1}|q|_{Q_{b_1}(P)},\\
    \norm{\bxi-\bPi_2^{0,k_2}\bxi}_{\bbM,P}\lesssim h_P^{s_2+1}|\bxi|_{\bV_2(P)},\quad \norm{\psi-\Pi_2^{0,k_2}\psi}_{Q_{b_2}(P)}\lesssim h_P^{s_2+1}|\psi|_{Q_{b_2}(P)}.
    \end{gather*}
\end{lemma}
\begin{lemma}\label{interpolation-estimates}
Given $\bv\in (\bH^{s_1+1}(P)\cap\bV_1(P), |\cdot|_{1,\bV_1(P)})$ and $\bxi \in (\bH^{s_2+1}(P)\cap \bV_2(P),|\cdot|_{1,\bV_2(P)})$. The Fortin interpolation operators $\bPi_1^{F,k_1}$ and $\bPi_2^{F,k_2}$ satisfy
\begin{align*}
    \norm{\bv-\bPi_1^{F,k_1}\bv}_{\bV_1(P)}\lesssim h_P^{s_1}|\bv|_{1,\bV_1(P)}, \quad \norm{\bxi - \bPi_2^{F,k_2}\bxi}_{\bbM,P} \lesssim h_P^{s_2+1}|\bxi|_{1,\bV_2(P)}.
\end{align*}
\end{lemma}
\paragraph{The virtual element formulation for the stress-assisted diffusion problem.} 
The discrete formulation for the fully-coupled problem reads: For given $\fb\in \bL^2(\Omega)$, $g\in L^2(\Omega)$, 
and $\varphi_D\in \text{H}^{1/2}(\Gamma_D)$, find 
$(\bu_h,p_h,\bzeta_h,\varphi_h) \in \bV_1^{h,k_1}\times Q_1^{h,k_1} \times \bV_2^{h,k_2} \times Q_2^{h,k_2}$ such that 
\begin{subequations}\label{eq:weak-discrete}
\begin{gather}
     \sum_{P\in \mathcal{T}^h} \biggl[2\mu \int_P \beps(\overline{\bu}_h):\beps(\overline{\bv}_h) + S_1^P(\bu_h-\overline{\bu}_h,\bv_h-\overline{\bv}_h)
     - \int_\Omega p_h\vdiv\bv_h \biggr] \notag \\ = \sum_{P\in \mathcal{T}^h} \int_P \overline{\fb}\cdot \bv_h, \; \forall \bv_h \in \bV_1^{h,k_1},\\
     \sum_{P\in \mathcal{T}^h} \left[ - \int_P q_h\vdiv\bu_h - \lambda^{-1} \int_P p_h q_h \right] = \sum_{P\in \mathcal{T}^h} \lambda^{-1}\int_P  \ell(\varphi_h)q_h, \; \forall q_h \in Q_1^{h,k_1},\\
     \sum_{P\in \mathcal{T}^h} \biggl[ \int_P \bbM(\overline{\bu}_h,p_h)^{-1} \overline{\bzeta}_h \cdot \overline{\bxi}_h + S_2^{\overline{\bu}_h,p_h,P}(\bzeta_h-\overline{\bzeta}_h,\bxi_h-\overline{\bzeta}_h) 
     + \int_P \varphi_h \vdiv \bxi_h \biggr] \notag \\
     = \sum_{P\in \mathcal{T}^h} \langle \varphi_D, \bxi_h\cdot\bn\rangle_{\partial P \cap \Gamma_D}, \; \forall \bxi_h \in \bV_2^{h,k_2},\\ 
     \sum_{P\in \mathcal{T}^h} \left[ \int_P \psi_h \vdiv\bzeta_h - \theta \int_\Omega \varphi_h\psi_h \right]  
    = - \sum_{P\in \mathcal{T}^h} \int_P g\psi_h, \forall \psi_h \in Q_2^{h,k_2},
\end{gather}
\end{subequations}
with $\overline{\bu}_h := \bPi_1^{\beps,k_1}\bu_h$, $\overline{\bv}_h := \bPi_1^{\beps,k_1}\bv_h$, $\overline{\fb}:=\bPi_1^{0,k_1-2}\fb$, $\overline{\bzeta}_h:=\bPi_2^{0,k_2} \bzeta_h$, $\overline{\bxi}_h:=\bPi_2^{0,k_2} \bxi_h$. The stabilisation terms are any symmetric and positive definite bilinear forms such that in the kernel of $\bPi_1^{\beps,k_1}$ (resp. $\bPi_2^{0,k_2}$) we have
$\norm{\bv_h}_{\bV_1} \lesssim S_{1}^{E}(\bv_h,\bv_h)\lesssim \norm{\bv_h}_{\bV_1(P)}$, $\norm{\bxi_h}_{\bbM,\Omega}\lesssim S_{2}^{\overline{\bu}_h,{p}_h,E}(\bxi_h,\bxi_h)\lesssim \norm{\bxi_h}_{\bbM,\Omega}$. 

Similarly to the continuous case, the continuous dependence on data for \eqref{eq:weak-discrete} follows as an extension of \cite[Section 4]{khot2024} (see Theorem~\ref{well-posedness}), this is a consequence of the stabilisation operators, and the Fortin-like interpolation operators that lead to the discrete inf-sup condition. Note the associated constants $\overline{C}_1,\overline{C}_2$ that appear in the discrete well-posedness (independent of the physical parameters) do not need to coincide with $C_1$, and $C_2$. We finalise by recalling the convergence result of the VE scheme.
\begin{theorem}\label{convergence-rates}
Under the assumptions of Theorem~\ref{well-posedness}. Given $(\bu,{p},\bzeta,\varphi)\in (\bH^{s_1+1}(\Omega)\cap \bV_1)\times (\text{H}^{s_1}(\Omega)\cap Q_{b_1}) \times (\bH^{s_2}\cap \bV_2) \times (\text{H}^{s_2}\cap Q_{b_2})$, $(\bu_h,{p}_h,\bzeta_h,\varphi_h)\in \bV_1^{h,k_1}\times Q_1^{h,k_1}\times \bV_2^{h,k_2}\times Q_2^{h,k_2}$ be the respective solutions of the continuous and discrete problems, with the data satisfying $\fb\in \bH^{s_1-1}\cap \bQ_{b_1}$ and $g\in H^{s_2}(\Omega)\cap Q_{b_2}$. If $\overline{C}_1 \sqrt{M} L_\ell + \overline{C}_2^2 \sqrt{M^3} L_\bbM\sqrt{2\mu}   (\norm{\varphi_D}_{1/2,\Gamma_D} + \norm{g}_{0,\Omega}) < 1/2.$ Then, the total error $\overline{\textnormal{e}}_h:=\norm{(\bu-\bu_h,{p}-{p}_h, \bzeta-\bzeta_h,\varphi-\varphi_h)}_{\bV_1\times Q_{1} \times \bV_2\times Q_2}$ decays with the following rate for $s:= \min \left\{s_1,s_2\right\}$
    \begin{align*}\label{convergence-rate}
         \overline{\textnormal{e}}_h &\lesssim h^{ s} (|\fb|_{s_1-1,\bQ_{b_1}} + |\bu|_{s_1+1,\bV_1} + |{p}|_{s_1,Q_{b_1}} + |g|_{s_2,Q_{b_2}} + |\bzeta|_{s_2,\bV_2}+|\varphi|_{s_2,Q_{b_2}}).
    \end{align*}
\end{theorem}

\section{Numerical results}\label{sec:results}
\begin{figure}[ht!]\label{fig:example-1}
    \centering
    \includegraphics[width=0.875\textwidth]{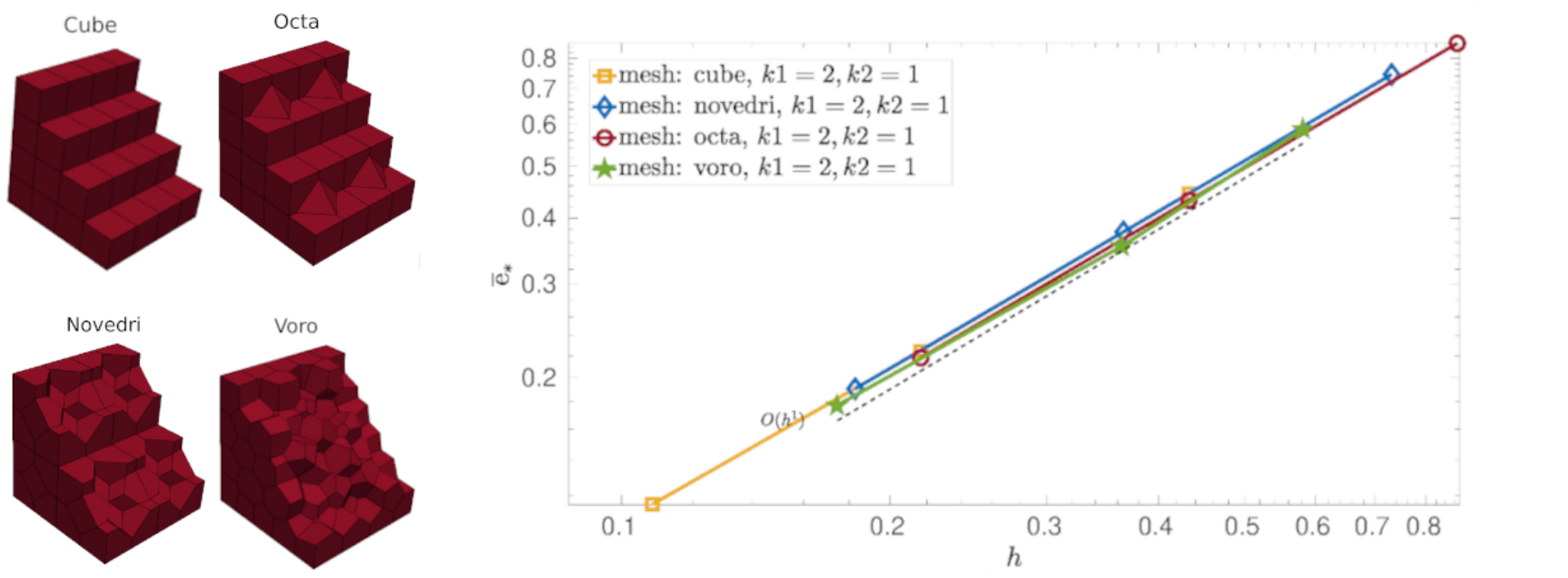}
    \caption{Example 1. An illustration of the distinct meshes used (left), and the converge plot under uniform refinement (right).}
\end{figure}

\begin{figure}[ht!]\label{fig:example-2}
    \centering
    \includegraphics[width=0.24\textwidth]{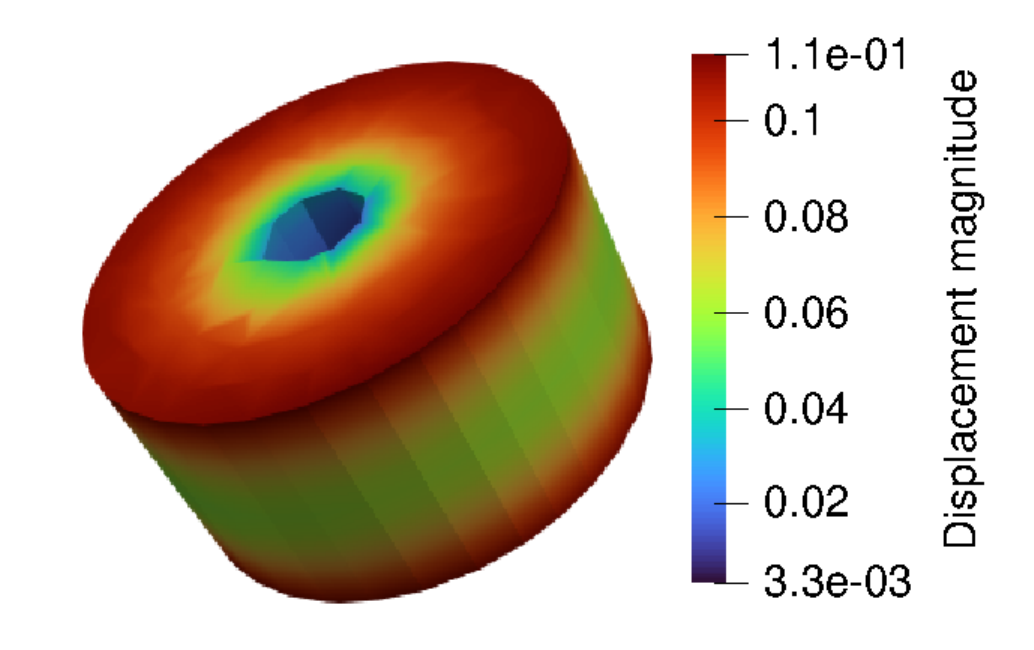}
    \includegraphics[width=0.24\textwidth]{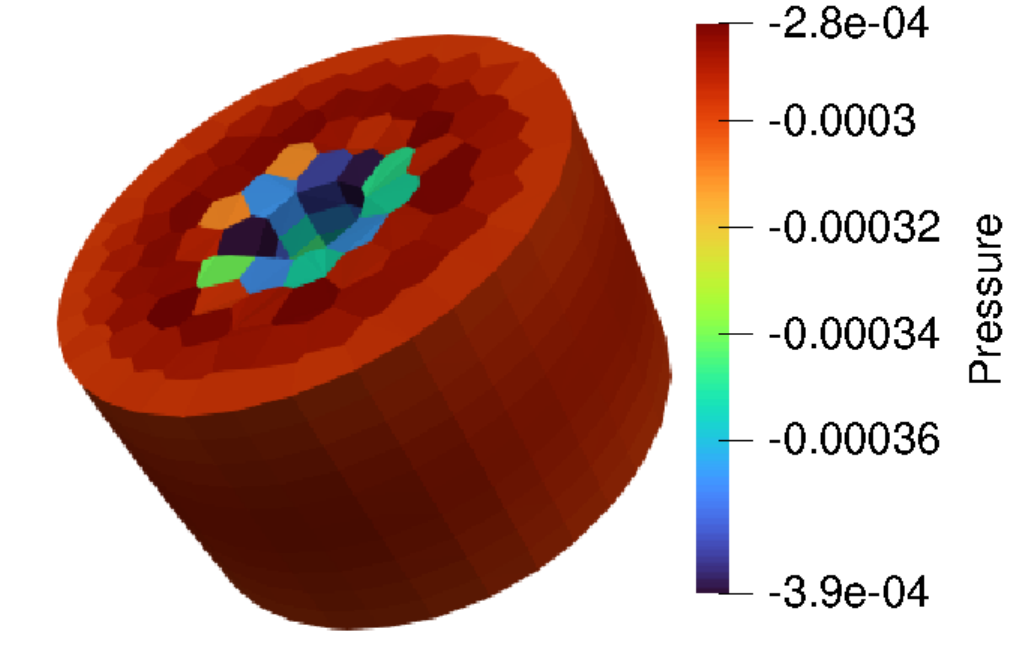}
    \includegraphics[width=0.24\textwidth]{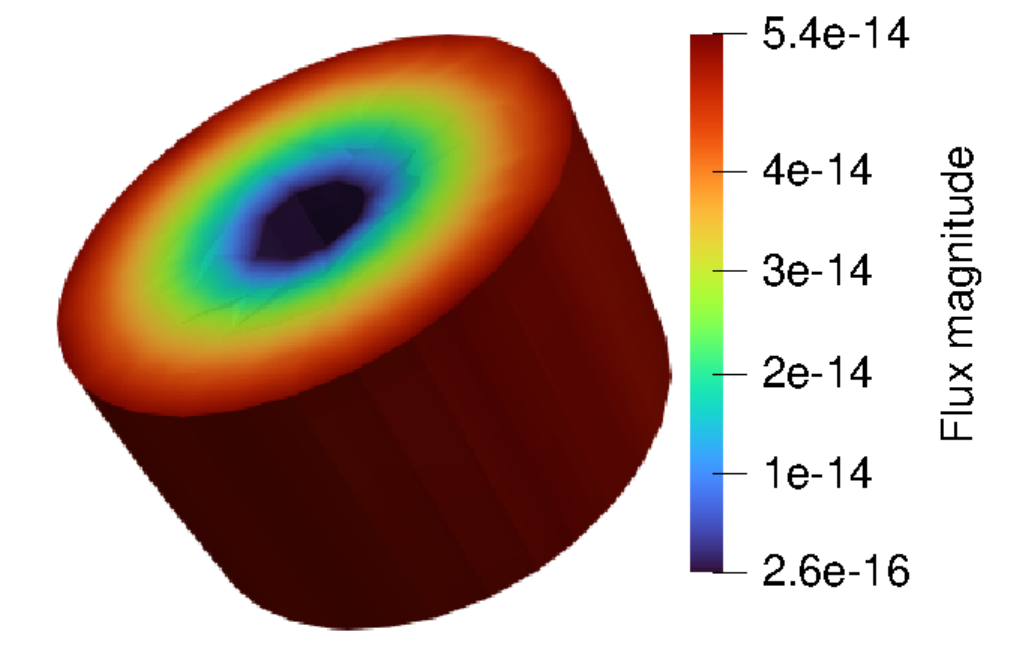} \includegraphics[width=0.24\textwidth]{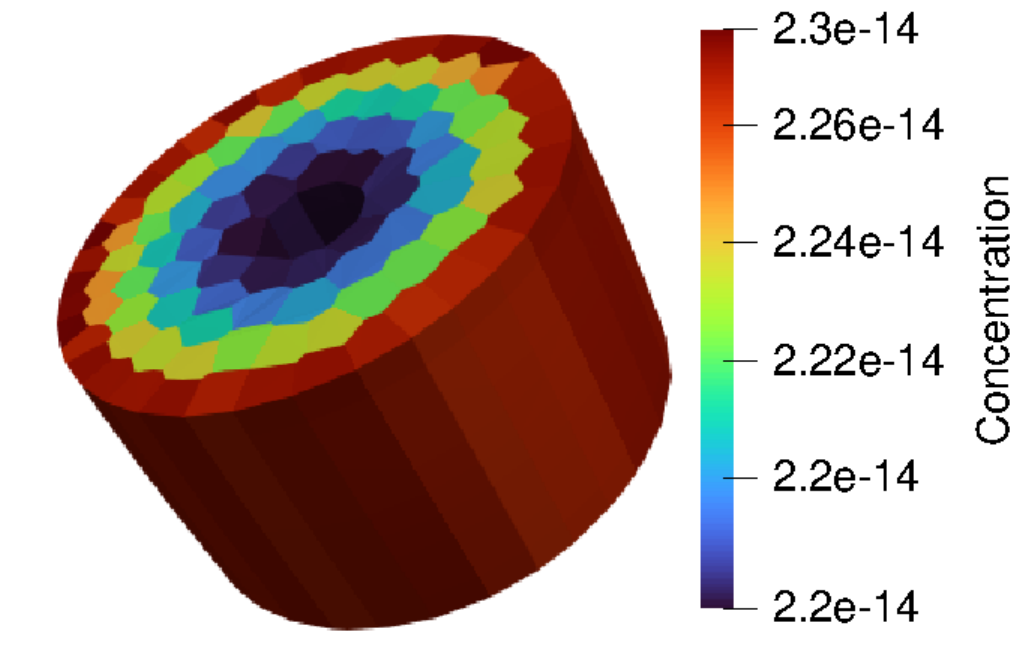}
    \includegraphics[width=0.24\textwidth]{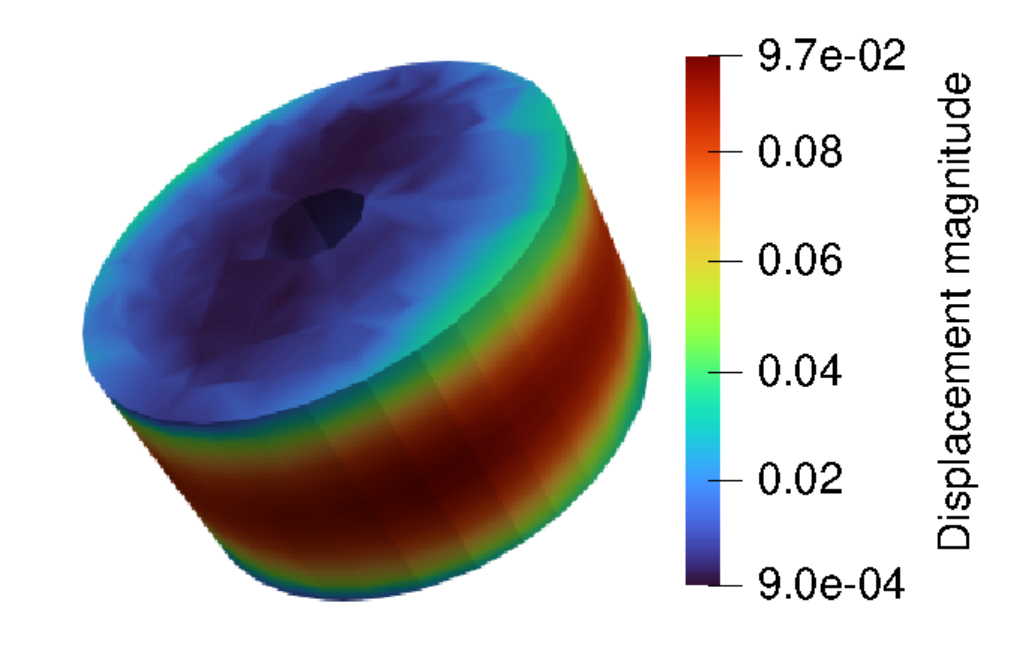}
    \includegraphics[width=0.24\textwidth]{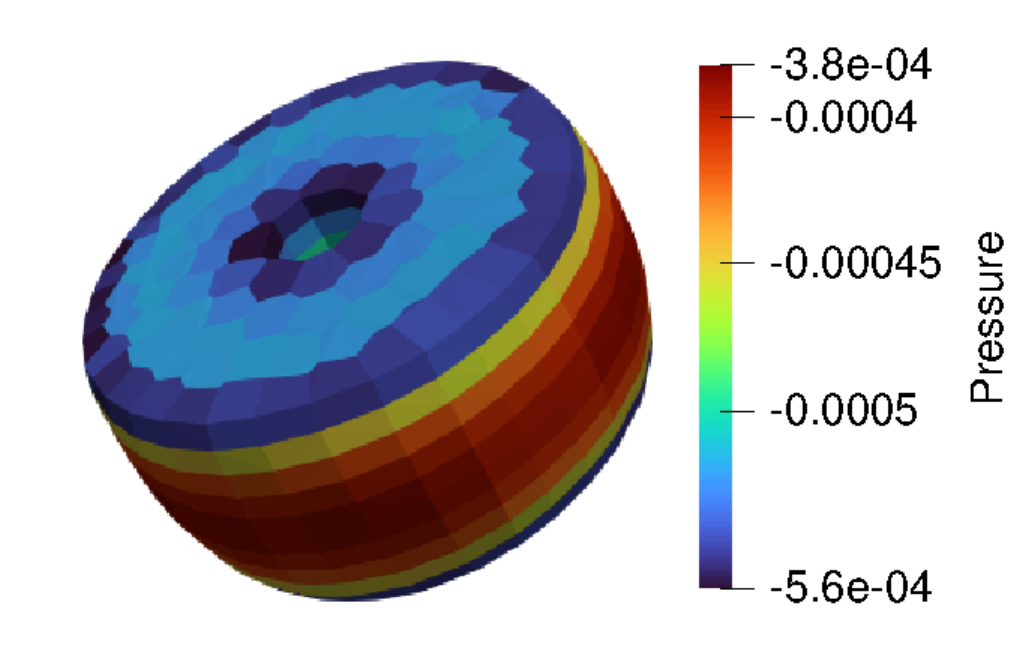}
    \includegraphics[width=0.24\textwidth]{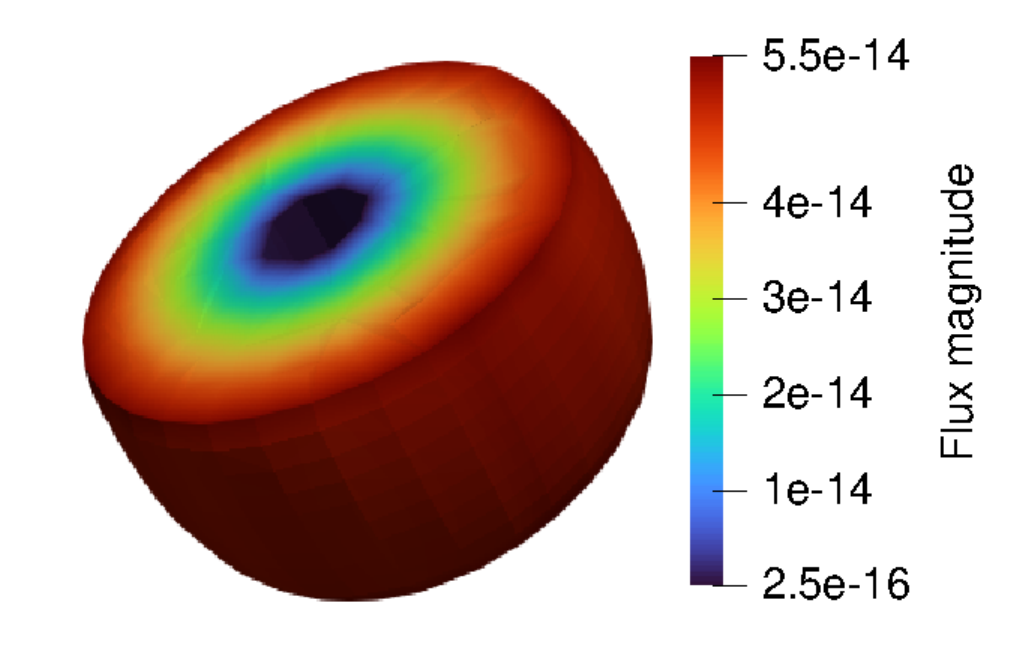} \includegraphics[width=0.24\textwidth]{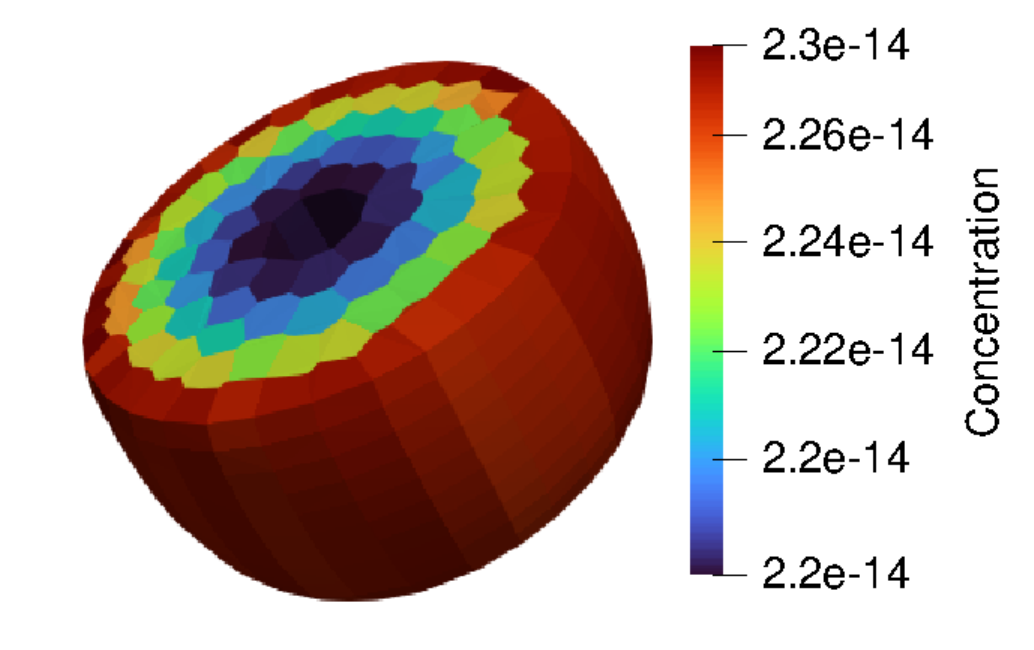}
    \caption{Example 2. Snapshots of the variables of interest in the reference configuration (left column), and deformed configuration (remaining columns) for the unclamped (first row), and clamped cases (second row). }
\end{figure}

 The simulations presented are implemented in VEM++ \cite{dassi2023vem++} for lowest case order ($k_1=2,\, k_2=1$). The fixed point algorithm based on the analysis in \cite{khot2024} is optimised to ensure that, in each iteration, only the blocks containing non-linearities are reconstructed, improving the performance of the code. We recall that the computational error is computed as $\overline{\textnormal{e}}_*:=\norm{(\bu-\overline{\bu}_h,{p}-{p}_h, \bzeta-\overline{\bzeta}_h,\varphi-\varphi_h)}_{\bV_1\times Q_{1} \times \bV_2\times Q_2}$, and the fixed point tolerance is set to $10^{-5}$.
\paragraph{Example 1.} Let $\Omega = (0,1)^3$ with $\Gamma_N = \{(x,y,z)\in \bbR^3 : x,y,z = 1\}$, and $\Gamma_D = \partial \Omega \setminus \Gamma_N$, we set the manufactured solutions as follows
$\bu(x,y,z) \hspace{-0.1cm} = \hspace{-0.1cm} \left( x^2 + x \cos(x) \sin(y), y^2 + x \cos(y) \sin(x), z^2 + x \cos(x) \cos(y) \right)/5$,
$\varphi(x,y,z) = \cos(\pi y) + \sin(\pi x) + x^2 + y^2 + z^2$,
$\ell(\varphi) = 1+\varphi^2/(1+\varphi^2)$, the non-linearities are given by
$\bbM(\beps(\bu),p) = 10^{-3}[ e^{-10^{-4}\tr((2\mu\beps(\bu)-p)\bbI)}]\bbI$,
with the adimensional parameters $\mu=10^2$, $\lambda=10^3$, $\theta=10^{-3}$, $M=2\times 10$. Figure~\ref{fig:example-1} confirms the linear convergence rate as predicted in Theorem~\ref{convergence-rates} for a variety of meshes, thanks to the flexibility provided by the VEM scheme proposed.
\paragraph{Example 2.} The lithiation process is configured in a perfored cylindrical particle of outer radious of $\qty{5}{\um}$, inner radious of $\qty{1}{\um}$, and height of $\qty{5}{\um}$. The boundary conditions are given by: Clamped and zero lithium fluxes on the inner circumference, the maximum lithium concentration $\hat{\Omega}=\qty{2.29d-14}{\mole/\mu \m^3}$ is fixed on the outer circumference, and a traction of $\qty{-2d-4}{\N/\mu \m^2}$ on the outer circumference. The Young's modulus is given by $E = \qty{1d-2}{\N/\mu\m^2}$ and the Poisson ratio is $\nu=0.3$. The diffusive source is zero and there is no body load force. The non-linear terms are given by $\bbM(\beps(\bu),p)= m_0(\bbI+m_0m_1 ((2\mu\beps(\bu)-p)\bbI)^2)$, and  $\ell(\varphi)=K_0 \varphi$, with $m_0=\qty{1d2}{\mu \m^2/\s}$, $m_1= \qty{1d3}{\mu \m^2/\N}$, $K_0=\tilde{\Omega}(2\mu+3\lambda)/3$ and the partial molar volume is $\tilde{\Omega}=\qty{3.497d12}{\mu \m^3/\mole}$. Finally, $\theta = M = 1$. In this experiment, we checked the behaviour of the solution when the top and bottom bases were unclamped and clamped (see Figure~\ref{fig:example-2}), the results coincide with the behaviour expected (see e.g. \cite{Taralov2015}), confirming the applicability of the model.

\paragraph{Acknowledgements.} This work has been partially supported by the Australian Research Council through the Future Fellowship Grant FT220100496. The author also thanks Prof. Ricardo Ruiz-Baier, Dr. Rekha Khot, and A/Prof. Franco Dassi for their constant support and expert guidance.

\small
\bibliography{plain}
\end{document}